# Quadratic Primes

N. A. Carella

**Abstract:** The subset of quadratic primes $\{ p = an^2 + bn + c : n \in \mathbb{Z} \}$ generated by some irreducible polynomial $f(x) = ax^2 + bx + c \in \mathbb{Z}[x]$ over the integers $\mathbb{Z}$ is widely believed to be an unbounded subset of prime numbers. This work provides the details of a possible proof for some quadratic polynomials $f(x) = x^2 + d$, $1 \leq d \leq 100$. In particular, it is shown that the cardinality of the simplest subset of quadratic primes $\{ p = n^2 + 1 : n \in \mathbb{Z} \}$ is infinite.

.

## 1. Introduction

Let $f(x) = ax^2 + bx + c \in \mathbb{Z}[x]$ be an irreducible polynomial over the integers $\mathbb{Z} = \{ ..., -3, -2, -1, 0, 1, 2, 3, ... \}$. The quadratic primes conjecture claims that certain Diophantine equation $y = ax^2 + bx + c$ has infinitely many prime solutions $y = p$ as the integers $n \in \mathbb{Z}$ varies over the set of integers. More generally, the Bouniakowsky conjecture claims that for an irreducible $f(x) \in \mathbb{Z}[x]$ over the integers of fixed divisor $\text{div}(f) = 1$, and degree $\deg(f) \geq 2$, the Diophantine equation $y = f(x)$ has infinitely many prime solutions $y = p$ as the integer $n \in \mathbb{Z}$ varies over the set of integers. The fixed divisor $\text{div}(f) = \gcd(f(\mathbb{Z}))$ of an irreducible polynomial $f(x) \in \mathbb{Z}[x]$ over the integers is the greatest common divisors of its image $f(\mathbb{Z}) = \{ f(n) : n \in \mathbb{Z} \}$ over the integers. The fixed divisor $\text{div}(f) = 1$ if the congruence $f(n) \equiv 0 \bmod p$ has $w(p) < p$ solutions for all prime numbers $p \leq \deg(f)$, see [FI10, p. 395]. Detailed discussions of the quadratic primes conjecture appear in [RN96, p. 387], [LG10, p. 17], [FI10, p. 395], [NW00, p. 405], [PJ09, p. 33], [IH78], [DI82], [LR12], [HW08], and related topics in [BZ07], [CC00], [GM00], [MK09], et alii.

This note provides the details of a possible proof for some quadratic polynomials $f(x) = x^2 + d$, $1 \leq d \leq 100$. In particular, the cardinality of the simplest subset of quadratic primes $\{ p = n^2 + 1 : n \in \mathbb{Z} \}$ is infinite. The techniques employed here are much simpler than the standard sieve methods employed in [IH78], [DI82], and [LR12], and a few other authors. This analysis is based on a simple weighted sieve. Essentially, it is a synthesis of those techniques used in [HB10, p. 1-4], and by other authors.

***Theorem* 1.** Let $f(x) = x^2 + d \in \mathbb{Z}[x]$ be an irreducible polynomial over the integers of fixed divisor $\text{div}(f) = 1$, and $1 \leq d \leq 100$. Then, the Diophantine equation $p = n^2 + d$ has infinitely many primes solutions $y = p$ as $n \in \mathbb{Z}$ varies over the integers.

***Proof***: Without loss in generality, let $f(x) = x^2 + 1$. Since the congruence $n^2 + 1 \equiv 0 \bmod p$ has less than $p$ solutions for any prime $p \leq \deg(f) = 2$, the fixed divisor of the polynomial is $\text{div}(f) = 1$. Now, select an appropriate weighted finite sum over the integers as observed in (16):

$$\sum_{n^2+1 \leq x} \frac{\Lambda(n^2+1)}{n\sqrt{\log n}} = -\sum_{n^2+1 \leq x} \frac{1}{n\sqrt{\log n}} \sum_{d \mid n^2+1} \mu(d) \log d$$

$$= -\sum_{d \leq x} \mu(d) \log d \sum_{\substack{n^2+1 \leq x, \\ n^2+1 \equiv 0 \bmod d}} \frac{1}{n\sqrt{\log n}}. \quad (1)$$

where $\Lambda(n) = -\sum_{d \mid n} \mu(d) \log d$, and $\gcd(d, n) = 1$. This follows from Lemma 2, and inverting the order of summation. Other examples of inverting the summation are given in [MV07, p. 35], [RH94, p. 27], [RM08, p. 216], [SN83, p. 83], [TM95, p. 36], and other.

Since the small moduli $d \geq 1$ contribute the bulk of the main term of the finite sum (1), and the large moduli have insignificant contribution, the last finite sum is broken up into two finite sums according to $d \leq x^{1/2-\epsilon}$ or $x^{1/2-\epsilon} < d \leq x$, where $\epsilon > 0$ is an arbitrarily small number. Specifically, the dyadic decomposition has the form

$$\sum_{n^2+1 \leq x} \frac{\Lambda(n^2+1)}{n\sqrt{\log n}} = -\sum_{d \leq x^{1/2-\epsilon}} \mu(d) \log d \sum_{\substack{n^2+1 \leq x, \\ n^2+1 \equiv 0 \bmod d}} \frac{1}{n\sqrt{\log n}} - \sum_{x^{1/2-\epsilon} < d \leq x} \mu(d) \log d \sum_{\substack{n^2+1 \leq x, \\ n^2+1 \equiv 0 \bmod d}} \frac{1}{n\sqrt{\log n}}. \quad (2)$$

**Small Moduli $d \leq x^{1/2-\epsilon}$, $\epsilon > 0$:** Applying Lemmas 4 and 7 yield

$$-\sum_{d \leq x^{1/2-\epsilon}} \mu(d) \log d \sum_{\substack{n^2+1 \leq x, \\ n^2+1 \equiv 0 \bmod q}} \frac{1}{n\sqrt{\log n}} \geq -\sum_{d \leq x^{1/2-\epsilon}} \mu(d) \log d \left( c_0 \frac{\rho(d)}{d} \sqrt{\log x} \left( 1 + O\left(\frac{1}{x^\epsilon \log x}\right) \right) \right)$$

$$\geq -c_0 \sqrt{\log x} \sum_{d \leq x^{1/2-\epsilon}} \frac{\mu(d) \rho(d) \log d}{d} + O\left( \frac{1}{x^\epsilon \sqrt{\log x}} \sum_{d \leq x^{1/2-\epsilon}} \frac{\rho(d) \log d}{d} \right) \quad (3)$$

$$\geq c_0 \sqrt{\log x} \left( c_1 + O\left(e^{-c\sqrt{\log x}}\right) \right) + O\left( \frac{\log^2 x}{x^{\epsilon/2}} \right),$$

where $\rho(q) = \#\{ n \leq x^{1/2} : n^2 + 1 \equiv 0 \bmod q \} = O(x^{\epsilon/2})$, see (11), and $c_0 = \sqrt{2} - \sqrt{2(1-\epsilon)} > 0$, $c_1 > 0$, $c > 0$ are constants.

**Large Moduli $d > x^{1/2-\epsilon}$, $\epsilon > 0$:** Applying Lemmas 5 and 6 yield

$$-\sum_{x^{1/2-\epsilon} < d \leq x} \mu(d) \log d \sum_{\substack{n^2+1 \leq x, \\ n^2+1 \equiv 0 \bmod q}} \frac{1}{n\sqrt{\log n}} = O\left(e^{-c\sqrt{\log x}} \log^4 x\right) + O\left(\frac{\log^4 x}{x^{(\log\log x)(\log\log\log x)/\log x}}\right). \tag{4}$$

Combining these expressions into (2) yield

$$\sum_{n^2+1 \leq x} \frac{\Lambda(n^2+1)}{n\sqrt{\log n}} \gg \sqrt{\log x}\left(1 + O\left(\frac{\log^{7/2} x}{x^{(\log\log x)(\log\log\log x)/\log x}}\right)\right). \tag{5}$$

Since the contribution by the subset of prime powers $n^2 + 1 = p^v \leq x$, $v \geq 2$, is zero, that is,

$$\sum_{n^2+1=p^v \leq x,\, v \geq 2} \frac{\Lambda(n^2+1)}{n\sqrt{\log n}} = 0, \tag{6}$$

see Lemmas 12 and 13 in Section 7, it follows that the cardinality of the subset of primes $\{n^2 + 1 = p \leq x, n \in \mathbb{Z}\}$ is unbounded as $x \longrightarrow \infty$. ∎

Some irreducible polynomials $f(x) = ax^2 + bx + c \in \mathbb{Z}[x]$ of fixed divisor $\text{div}(f) = 1$ can be transformed into an equivalent case of the form $g(x) = x^2 + d$ by means of algebraic manipulations. The equivalent problem is then handled as the case $f(x) = x^2 + 1$ presented above. Moreover, a lower bound of the asymptotic formula is provided in Section 8.

The topics of primes in quadratic, cubic, quartic arithmetic progressions, the Bouniakowsky conjecture, and in general the Hypothesis H, and the Bateman-Horn conjecture, [RN96], are rich areas of research involving class fields theory and analytic number theory. The linear case, Dirichlet Theorem for primes in arithmetic progressions, is proved in [CS09, Chapter 2] from the point of view of class fields theory. The quadratic case, Theorem 1 here, seems to have another proof in term of Hecke $L$-functions over the Gaussian quadratic field $\mathbb{Q}(\sqrt{-1})$ quite similar to Dirichlet Theorem's proof in term of $L$-functions over the rational field $\mathbb{Q}$.

As an application, it can be shown that there are irreducible quadratic polynomials $f(x) = ax^2 + bx + c \in \mathbb{Z}[x]$ of fixed divisor $\text{div}(f) \neq 1$ such that $f(n) = P_2(n)$ is the product of two primes for infinitely many values. For example, take the irreducible polynomial $f(x) = x(x+1) + 2$ of fixed divisor $\text{div}(f) = 2$. By Theorem 1, $f(n)/\text{div}(f) = p$ is prime for infinitely many values. Ergo, $f(n) = 2^2 p$ for infinitely many values $n \in \mathbb{Z}$. The topic of almost primes is studied in [IH78], [DI82], and [LR12], using sieve methods.

The elementary underpinning of Theorem 1 is assembled in Sections 2 to 8. In Lemma 14, a lower bound of the correspond-ing primes counting function is computed. The proofs of all these lemmas use elementary methods.

## 2. Elementary Foundation

The basic definitions of several number theoretical functions, and a handful of Lemmas are recorded here.

### 2.1 Formulae for the Mobius and vonMangoldt Functions

Let $n \in \mathbb{N} = \{0, 1, 2, 3, \ldots\}$ be an integer. The Mobius function is defined by

$$\mu(n) = \begin{cases} (-1)^t & \text{if } n = p_1^{v_1} \cdot p_2^{v_2} \cdots p_t^{v_t}, \text{ if } v_1 = \cdots = v_t = 1, \\ 0 & \text{if } n = p_1^{v_1} \cdot p_2^{v_2} \cdots p_t^{v_t}, \text{ if some } v_i \neq 1. \end{cases} \quad (7)$$

The subset of squarefree numbers $\{n \in \mathbb{N} : \mu(n) \pm 1\} = \{n = p_1 \cdot p_2 \cdots p_t : p_i \text{ prime}\}$ is the support of the Mobius function $\mu : \mathbb{N} \longrightarrow \{-1, 0, 1\}$. Further, the vonMangoldt function is defined by

$$\Lambda(n) = \begin{cases} \log p & \text{if } n = p^k, k \geq 1, \\ 0 & \text{if } n \neq p^k, k \geq 1. \end{cases} \quad (8)$$

The subset of prime powers $\{n \in \mathbb{N} : \Lambda(n) \neq 0\} = \{n \in \mathbb{N} : n = p^k, k \geq 1\}$ is the support of the vonMangoldt function $\Lambda : \mathbb{N} \longrightarrow \mathbb{R}$.

***Lemma 2.*** Let $n \geq 1$ be an integer, and let $\Lambda$ be the vonMangoldt function. Then

$$\Lambda(n) = -\sum_{d \mid n} \mu(d) \log d \quad (9)$$

***Proof***: Use Mobius inversion formula

$$f(n) = \sum_{d \mid n} g(d) \iff g(n) = \sum_{d \mid n} \mu(d) f(n/d) \quad (10)$$

on the identity $\log n = \sum_{d \mid n} \Lambda(d)$ to confirm this claim. ∎

Extensive details for other identities and approximations of the vonMangoldt are discussed in [GP05], [FI10], [HN07, p. 27], et alii.

### 2.2 Quadratic Equations Over Arithmetic Progressions

Let $1 \leq a < q$ be integers, $\gcd(a, q) = 1$. An element $a \in \{0, 1, 2, 3, \ldots q-1\} \cong \mathbb{Z}/q\mathbb{Z}$ is called a quadratic residue if the congruence $z^2 \equiv a \bmod q$ has a solution. Otherwise, $a \geq 1$ is a quadratic nonresidue. For a quadratic residue $a$ modulo $q$, with $\gcd(a, q) = 1$, the congruence $z^2 \equiv a \bmod q$ has $2^W \geq 2$ solutions, where $W = \omega(q) + r$, $r = 0, 1, 2$, according as $4 \nmid q$, or $4 \| q$ or $8 \mid q$, see (11) and [LV56, p. 65]. The function $\omega(q)$ tallies the number of prime divisors of $q$, cf [CC07]

**Lemma 3.** (Fermat) Let $p \geq 3$ be a prime number. Then

(i) The integer $-1$ is quadratic residue modulo $p$ if and only if $p = u^2 + v^2$, $u, v \in \mathbb{Z}$.

(ii) The integer $-1$ is quadratic nonresidue modulo $p$ if and only if $p \neq u^2 + v^2$, $u, v \in \mathbb{Z}$.

The quadratic residuacity of numbers and related concepts are explicated in almost every textbook in elementary number theory, [HW08], [SN83], [LV56], [RH94], et alii.

The previous information quickly leads to a formula for the multiplicative function $\rho(q) = \#\{n \leq x^{1/2} : n^2 + 1 \equiv 0 \bmod q\}$. For any integer $q \geq 2$ this is defined by

$$\prod_{p \mid q} \rho(p), \quad \text{where} \quad \rho(p) = \begin{cases} 1 & \text{if } p = 2, \\ 2 & \text{if } p = u^2 + v^2, \\ 0 & \text{if } p \neq u^2 + v^2 \text{ or } p = 2^k, k \geq 2. \end{cases} \tag{11}$$

These concepts come into play in the analysis of powers sums, and other finite sums over quadratic arithmetic progressions.

## 3. Finite Sums Over Small Moduli

The small moduli $q = O(\sqrt{\log x})$, as demonstrated below, have the largest effect and contribute the bulk of the main term in the finite sum considered here.

**Lemma 4.** Let $x \in \mathbb{R}$ be a large number, and let $q < x$ be an integer. Then,

(i) $\displaystyle\sum_{\substack{n^2+1 \leq x, \\ n^2+1 \equiv 0 \bmod q}} \frac{1}{n\sqrt{\log n}} = \frac{2\rho(q)}{q}\left(\sqrt{\log(x^{1/2} - q f(q))} - \sqrt{\log r(q)}\right),$ if $q < x$.

(ii) $\displaystyle\sum_{\substack{n^2+1 \leq x, \\ n^2+1 \equiv 0 \bmod q}} \frac{1}{n\sqrt{\log n}} \geq c_0 \frac{\rho(q)}{q} \sqrt{\log x}\left(1 + O\left(\frac{1}{x^\epsilon \log x}\right)\right),$ if $q \leq x^{1/2-\epsilon}$, $\epsilon > 0$.

(12)

where $r(q) \geq 1$ is a solution of the quadratic congruence $n^2 + 1 \equiv 0 \bmod q$, the term $\rho(q) = \#\{n \leq x^{1/2} : n^2 + 1 \equiv 0 \bmod q\}$ is the number of solutions, the term $f(q) = \{(x^{1/2} - r(q))/q\}$ is the fractional part function, and $c_0 = \sqrt{2} - \sqrt{2(1-\epsilon)}$ is a constant.

**Proof**: (i) Fix a solution $r = r(q) \geq 1$ of the quadratic congruence $n^2 + 1 \equiv 0 \bmod q$ such that $r^2 + 1 \leq x$. Then, each integer $n = qm + r \leq x^{1/2}$ is also a solution. Thus, the finite sum can be rewritten as

$$\sum_{\substack{n^2+1 \leq x, \\ n^2+1 \equiv 0 \bmod q}} \frac{1}{n\sqrt{\log n}} = \sum_{m \leq V} \frac{1}{(qm+r)\sqrt{\log(qm+r)}},$$

(13)

where $V = \left(x^{1/2} - r\right)/q - \{(x^{1/2} - r)/q\}$, and $r = r(q) < q < x$. Evaluating the integral representation yields

$$\sum_{\substack{n^2+1 \leq x, \\ n^2+1 \equiv 0 \bmod q}} \frac{1}{n\sqrt{\log n}} = \int_0^V \frac{1}{(qt+r)\sqrt{\log(qt+r)}} \, dt$$

$$= \frac{2}{q}\sqrt{\log(qt+r)}\Big|_0^V = \frac{2}{q}\left(\sqrt{\log(x^{1/2} - qf(q))} - \sqrt{\log r(q)}\right),$$

(14)

where $f(q) = \{(x^{1/2} - r)/q\}$ is the fractional part function.

For (ii), let $q \leq x^{1/2-\epsilon}$, where $\epsilon > 0$ is an arbitrarily small number, and $r(q) < q$, then

$$\sqrt{\log(x^{1/2} - qf(q))} - \sqrt{\log r(q)} \geq \sqrt{\log(x^{1/2} - x^{1/2-\epsilon})} - \sqrt{\log x^{1/2-\epsilon}}$$

$$\geq \sqrt{\log x^{1/2}\left(1 + \frac{\log(1 - x^{-\epsilon})}{\log x^{1/2}}\right)} - \sqrt{\log x^{1/2-\epsilon}}$$

$$\geq \left(\sqrt{1/2} - \sqrt{1/2 - \epsilon}\right)\sqrt{\log x} + O\left(\frac{1}{x^\epsilon \sqrt{\log x}}\right)$$

$$> 0.$$

(15)

Substituting this back, and simplifying it, provide a lower bound of the main term. ∎

Various other finite sums having unbounded partial sums appear to be suitable for this analysis. For example,

$$\sum_{n^2+1 \le x} \frac{\Lambda(n^2+1)}{n(\log n)^{1-\alpha}}, \quad \alpha > 0. \tag{16}$$

## 4. Finite Sums Over Large Moduli

This Section is concerned with effective estimates of certain finite sums over large moduli. The finite sum over the subset of large moduli, which have small number of prime divisors $\omega(q) \le \log\log x$, is covered in Subsection 4.1. The other finite sum over the subset of moduli, which have large number of prime divisors $\omega(q) > \log\log x$, is covered in Subsection 4.2.

### 4.1 Case $\omega(q) \le \log\log x$

Let $\omega(q)$ be the number of prime divisors of $q$. Each integer in the subset of large moduli $\{ q \le x : \omega(q) \le \log\log x \}$ has a small number of prime divisors. Accordingly, the quadratic congruence $z^2 \equiv a \bmod q$ has a small number of solutions $\rho(q) \le 2^{\omega(q)+2} \le \log x$. The subset of integers $\{ q \le x : \omega(q) \le \log\log x \}$ has density 1 in the set of nonnegative integers $\mathbb{N} = \{ 0, 1, 2, 3, \dots \}$.

**Lemma 5.** Let $x \in \mathbb{R}$ be a large number, and let $q > x^{1/2-\epsilon}$, $\epsilon > 0$, be an integer such that $\omega(q) \le \log\log x$. Then,

$$-\sum_{\substack{x^{1/2-\epsilon} < q \le x, \\ \omega(q) \le \log\log x}} \mu(q) \log q \sum_{\substack{n^2+1 \le x, \\ n^2+1 \equiv 0 \bmod q}} \frac{1}{n\sqrt{\log n}} \ll e^{-c\sqrt{\log x}} \log^2 x. \tag{17}$$

**Proof**: Let $r \ge 1$ be a root of $n^2 + 1 \equiv 0 \bmod q$. As the moduli $q$ are restricted to the range $x^{1/2-\epsilon} < q \le x$, and $n = qm + r \le x^{1/2}$, these data imply that $x^{1/2-\epsilon} m + r \le n = qm + r \le x^{1/2}$. Hence, the index $m$ is restricted to the range $0 \le m \le x^{\epsilon}$.

Apply Lemma 2, see also (1) for finer details, and rearrange the finite sum as

$$-\sum_{\substack{x^{1/2-\epsilon} < q \le x, \\ \omega(q) \le \log\log x}} \mu(q) \log q \sum_{\substack{n^2+1 \le x, \\ n^2+1 \equiv 0 \bmod q}} \frac{1}{n\sqrt{\log n}}$$

$$= -\sum_{\substack{x^{1/2-\epsilon} < q \le x, \\ \omega(q) \le \log\log x}} \mu(q) \rho(q) \log q \sum_{m \le x^\epsilon} \frac{1}{(qm+r)\sqrt{\log(qm+r)}},$$

(18)

where $\rho(q) = \#\{ n \le x^{1/2} : n^2 + 1 \equiv 0 \bmod q \} \le \log x$. Use an integral representation to estimate the inner finite sum as follows:

$$\sum_{m \leq x^\epsilon} \frac{1}{(qm+r)\sqrt{\log(qm+r)}}$$

$$= \int_0^{x^\epsilon} \frac{1}{(qt+r)\sqrt{\log(qt+r)}} dt = \frac{2}{q}\left(\sqrt{\log(q \log x + r)} - \sqrt{\log r}\right),$$

(19)

where $1 \leq r < q < x$. Substituting these information, using summation by part, and applying Lemma 7, return

$$-\sum_{\substack{x^{1/2-\epsilon} < q \leq x, \\ \omega(q) \leq \log\log x}} \mu(q) \rho(q) \log q \sum_{\substack{n^2+1 \leq x, \\ n^2+1 \equiv 0 \bmod q}} \frac{1}{n\sqrt{\log n}}$$

$$= -\sum_{\substack{x^{1/2-\epsilon} < q \leq x, \\ \omega(q) \leq \log\log x}} \mu(q) \rho(q) \log q \left(\frac{2}{q}\left(\sqrt{\log(q x^\epsilon + r)} - \sqrt{\log r}\right)\right)$$

$$= -2 \sum_{\substack{x^{1/2-\epsilon} < q \leq x, \\ \omega(q) \leq \log\log x}} \frac{\mu(d) \rho(q) \log q}{q} \sqrt{\log(q x^\epsilon + r)} + 2 \sum_{\substack{x^{1/2-\epsilon} < q \leq x, \\ \omega(q) \leq \log\log x}} \frac{\mu(d) \rho(q) \log q}{q} \sqrt{\log r}$$

$$\ll e^{-c\sqrt{\log x}} \log^4 x,$$

(20)

where $r = r(q)$, refer to Lemma 7-iii, this estimate uses $\rho(q) \ll \log x$. ∎

The reader can confer [MV07, p. 182], and [RM08, p. 318] for similar evaluations.

### 4.2 Case $\omega(q) > \log\log x$

Let $\omega(q)$ be the number of prime divisors of $q$. Each integer in the subset of large moduli $\{q \leq x : \omega(q) > \log\log x\}$ has a large number of prime divisors. Accordingly, the quadratic congruence $z^2 \equiv a \bmod q$ has a large number of solutions $\rho(q) \leq 2^{\omega(q)+2} \leq q^\delta$, where $\delta > 0$ is an arbitrarily small number.

The subset $\{q \leq x : \omega(q) > \log\log x\}$ of such highly composite integers has zero density in the set of integers $\mathbb{N} = \{0, 1, 2, 3, ...\}$. In fact, this is a very small subset of integers $\#\{q \leq x : \omega(q) > \log\log x\} = o(x)$. This topic is discussed in Section 6, and some related information are given in [AE44, p. 449].

**Lemma 6.** Let $x \in \mathbb{R}$ be a large number, let $q > x^{1/2-\epsilon}$, $\epsilon > 0$, be an integer, and let $\omega(q) > \log\log x$. Then,

$$-\sum_{\substack{x^{1/2-\epsilon}<q\leq x,\\ \omega(q)>\log\log x}} \mu(d)\log d \sum_{\substack{n^2+1\leq x,\\ n^2+1\equiv 0 \bmod q}} \frac{1}{n\sqrt{\log n}} \ll \frac{\log^4 x}{x^{(\log\log x)(\log\log\log x)/\log x}}.$$

*Proof*: Here the moduli $q$ are restricted to the range $x^{1/2-\epsilon} < q \leq x$, and $n = qm + r \leq x^{1/2}$. Hence, the index $m$ is restricted to the range $0 \leq m \leq x^\epsilon$. Let $r(q) \geq 1$ be a root of $z^2 + 1 \equiv 0 \bmod q$. Next rearrange the finite sum as

$$-\sum_{\substack{x^{1/2-\epsilon}<q\leq x,\\ \omega(q)>\log\log x}} \mu(d)\log d \sum_{\substack{n^2+1\leq x,\\ n^2+1\equiv 0 \bmod q}} \frac{1}{n\sqrt{\log n}}$$

$$= -\sum_{\substack{x^{1/2-\epsilon}<q\leq x,\\ \omega(q)>\log\log x}} \mu(d)\rho(q)\log d \sum_{m\leq x^\epsilon} \frac{1}{(qm+r)\sqrt{\log(qm+r)}}$$

$$= -\sum_{\substack{x^{1/2}\log^{-1} x<q\leq x,\\ \omega(q)>\log\log x}} \mu(d)\rho(q)\log d\left(\frac{2}{q}\left(\sqrt{\log(qx^\epsilon+r)} - \sqrt{\log r(q)}\right)\right)$$

$$\ll \sqrt{\log x^\epsilon} \sum_{\substack{x^{1/2-\epsilon}<q\leq x,\\ \omega(q)>\log\log x}} \frac{\rho(q)\log^2 q}{q},$$

(22)

refer to Lemma 4 for similar calculations. Applying Lemma 11 completes the proof. ∎

## 5. Finite Sums of Mobius Function

**Lemma 7.** Let $x \in \mathbb{R}$ be a large number, and let $f(n) = O(\log^B x)$, $B > 0$, be a function. Let $q \geq 1$, $a \geq 1$ be a pair of integers, $\gcd(a, q) = 1$. Then,

(i) $\displaystyle\sum_{n\leq x,\, n\equiv a \bmod q} \mu(d) = O\left(x e^{-c\sqrt{\log x}}\right),$

(ii) $\displaystyle\sum_{n\leq x,\, n\equiv a \bmod q} \frac{\mu(n)\log n}{n} = -c_0 + O\left(e^{-c\sqrt{\log x}}\right),$ (23)

(iii) $\displaystyle\sum_{n\leq x,\, n\equiv a \bmod q} \frac{\mu(n)\log n\, f(n)}{n} = O\left(x e^{-c\sqrt{\log x}} f(x)\log x\right),$

where $c_0 > 0$, $c > 0$ are constants.

Some of these finite sums are evaluated in [MV07, p. 182-185], [NW00, p. 351], and [RM08].

## 6. Highly Composite Integers

The statistics on subsets of integers with specified number of primes is of general interest in many area of mathematics. Some information concerning these integers are recorded in this Section.

**Lemma 8.** (Landau) Let $x \in \mathbb{R}$ be a large number, and let $\pi_k(x) = \#\{n \leq x : \omega(n) = k\}$ be the counting function for the number of integers with $k$-prime factors. Then,

$$\pi_k(x) = \frac{x(\log\log x)^{k-1}}{(k-1)!(\log x)} + o\left(\frac{x(\log\log x)^{k-1}}{(k-1)!(\log x)}\right). \tag{24}$$

The proof is usually done for the restricted range $k \ll \log\log x$, but it holds for any parameter $k \leq \log x$. Detailed proofs are given in [DF12, p. 159], [ND12], and [TM95].

**Lemma 9.** Let $x \in \mathbb{R}$ be a large number, and let $\pi_k(x) = \#\{n \leq x : \omega(n) = k\}$, be the k-prime factors integers counting function. Then,

$$\sum_{\log\log x \leq k \leq \log x} \pi_k(x) \ll x. \tag{25}$$

*Proof*: By Lemma 8, this finite sum has the upper bound

$$\sum_{\log\log x \leq k \leq \log x} \pi_k(x) \ll \sum_{k \leq \log x} \frac{x(\log\log x)^{k-1}}{(k-1)!(\log x)} \ll \frac{x}{\log x} \sum_{k \geq 1} \frac{(\log\log x)^{k-1}}{(k-1)!} \tag{26}$$

$= x.$ ∎

The function $W_0(t) = \#\{q \leq t : \log\log t \leq \omega(q) \leq \log t\}$ accounts for the cardinality of the subset of integers with at least $k \geq \log\log t$ prime factors. An estimate for the closely related finite sum involving the counting measure $\rho_k(q) = \#\{n \leq x^{1/2} : n^2 + 1 \equiv 0 \bmod q \text{ and } \omega(q) = k\}$ is stated here.

**Lemma 10.** Let $x \in \mathbb{R}$ be a large number, and let $\pi_k(x) = \#\{n \leq x : \omega(n) = k\}$, be the $k$-prime factors integers counting function. Then,

$$W(x) = \sum_{\log\log x \leq k \leq \log x} \rho_k(x)\pi_k(x) \ll x^{1-\frac{(\log\log x)(\log\log\log x)}{2\log x}}.$$
(27)

*Proof*: Express the following quantities in powers of $x \geq 1$:

$$(\log x)^k = x^{\frac{k\log\log x}{\log x}}, \qquad\qquad (\log\log x)^{k-1} = x^{\frac{(k-1)\log\log\log x}{\log x}},$$

$$(k-1)! \ll x^{\frac{k\log k - k + O(\log k)}{\log x}}, \qquad\qquad \rho(q) \leq 2^{\omega(q)+2} = 2^{k+2} \leq x^{\frac{k+2}{\log x}}.$$
(28)

Substituting the previous expressions into the finite sum upper bound of the inner finite sum, and using Lemma 8, return

$$\sum_{\log\log x \leq k \leq \log x} \rho_k(x)\pi_k(x) \ll \sum_{\log\log x \leq k \leq \log x} 2^{k+2}\frac{x(\log\log x)^{k-1}}{(k-1)!(\log x)}$$

$$\ll \sum_{\log\log x \leq k \leq \log x}\left(x^{\frac{k+2}{\log x}}\right)\left(x^{1+\frac{(k-1)\log\log\log x - \log\log x - k\log k + k + O(\log k)}{\log x}}\right)$$

$$\ll x \sum_{\log\log x \leq k \leq \log x} x^{\frac{(k-1)\log\log\log x - \log\log x - k\log k + 2k + 2 + O(\log k)}{\log x}}$$

$$\ll x \sum_{\log\log x \leq k \leq \log x} x^{\frac{-k\log k}{\log x}} \ll x\frac{\log x}{x^{\frac{(\log\log x)(\log\log\log x)}{\log x}}}.$$
(29)

These complete the estimate. ∎

**Lemma 11.** Let $x \in \mathbb{R}$ be a large number, let $q > x^{1/2}\log^{-1}x$ be an integer, and let $\omega(q) > \log\log x$. Then,

$$\sum \frac{\rho(q)\log^m q}{q} \ll x^{1-\frac{(\log\log x)(\log\log\log x)}{2\log x}}, \quad m \geq 1.$$
(30)

*Proof*: Let $\rho_k(q) = \#\{n \leq x^{1/2} : n^2 + 1 \equiv 0 \bmod q \text{ and } \omega(q) = k\}$, and observe that

$$\sum_{\substack{x^{1/2-\epsilon} < q \leq x, \\ \omega(q) > \log\log x}} \rho(q) \leq \sum_{\substack{1 < q \leq x, \\ k > \log\log x}} \rho_k(q) = \sum_{k > \log\log x} \rho_k(x)\pi_k(x) = W(t).$$
(31)

Now, the finite sum is estimated as follows:

$$\sum_{\substack{x^{1/2-\epsilon}<q\leq x,\\ \omega(q)>\log\log x}} \frac{\rho(q)\log^m q}{q} \leq \sum_{\substack{1<q\leq x,\\ \omega(q)>\log\log x}} \frac{\log^m q}{q} = \int_1^x \frac{\log^m t}{t} dW(t)$$

$$= \left.\frac{\log^m t}{t} W(t)\right|_1^x - \int_1^x \left(\frac{\log^m t}{t}\right)' W(t)\, dt \qquad (32)$$

$$\ll \frac{\log^m x}{x} W(t).$$

Apply Lemma 10 to complete the estimate. ∎

## 7. Small Distances Between Powers Of Integers

The Catalan conjecture claims that the sequence of integers powers 1, $2^2$, $2^3$, $3^2$, $2^4$, $5^2$, $3^3$, $2^5$, $6^2$, $7^2$, $2^6$, $3^4$, $10^2$, ... has only one pair of consecutive powers, namely, $2^3$ and $3^2$. This result, usually expressed by the Diophantine equation $x^p - y^q = 1$, $p, q \geq 2$, was proved a few years ago, detailed information is widely available in the literature. The proofs of several special cases were established long ago.

**Lemma 12.** (Lebesgue-Nagell) For any exponent $n \in \mathbb{N}$, the Diophantine equation $x^2 + 1 = y^n$ has no nonzero integers $x, y \in \mathbb{Z}$ solutions.

For even $n = 2m$, $m \geq 1$ the proof is trivial. And for odd $n = 2m + 1$, the equation has genus $g > 1$, so it has a finite number of solutions. An algebraic proof appears in various places in the literature, [SC08, p. 9]. The generalized Lebesgue-Nagell equation has been completely solved for a small range of parameters.

**Lemma 13.** ([BS06]) The Diophantine equation $x^2 + d = y^n$, with $n \geq 3$, and $1 \leq d \leq 100$, has at most eight integer solutions $x, y \in \mathbb{Z}$.

A complete table of solutions appears in [BS06]. A few of these equations have no solutions at all, for example, $x^2 + 1 = y^n$ and $x^2 + 3 = y^n$. On the other direction, $x^2 + 28 = y^n$ has the most solutions:

$(x, y, n) = (6, 4, 3),\ (22, 8, 3),\ (225, 37, 3),\ (2, 2, 5),\ (6, 2, 6),\ (10, 2, 6),\ (22, 2, 9),\ (362, 2, 17)$.

Both Lemmas 12 and 13 immediately imply that the quadratic arithmetic progressions $\{n^2 + d : n \in \mathbb{Z}\}$ have bounded numbers of prime powers, that is, $n^2 + d = p^v : n \in \mathbb{Z}, v \geq 2$. Specifically, the finite sums

$$\sum_{n^2+1=p^v \leq x,\, v \geq 2} \frac{\Lambda(n^2 + 1)}{n\sqrt{\log n}} = 0 \quad \text{and} \quad \sum_{n^2+d=p^v \leq x,\, v \geq 2} \frac{\Lambda(n^2 + d)}{n\sqrt{\log n}} = O(1) \tag{33}$$

for $1 \leq d \leq 100$. Employing standard analytical method, it can be shown that these finite sums are bounded by a constant. But the algebraic proofs of the Lebesgue-Nagell equation give exact answers.

The number of solutions of the more general case $n^2 + d = p^v$, $n \in \mathbb{Z}$, $v \geq 2$, and $d \neq 0$ constant, will be of interest in proof for primes of the form $n^2 + d = p : n \in \mathbb{Z}$.

## 8. Asymptotic Formula

Let $x \in \mathbb{R}$ be a large number, and let $\pi_f(x) = \#\{f(n) \leq x : f(n) = p \text{ is prime}\}$ be the corresponding counting function of the prime numbers defined by the irreducible polynomial $f(x) \in \mathbb{Z}[x]$. For the specific case $f(x) = x^2 + 1$, the expected asymptotic formula has the form

$$\pi_f(x) = \prod_{p \geq 3}\left(1 - \frac{1}{p-1}\left(\frac{-1}{p}\right)\right)\frac{\sqrt{x}}{\log x} + o\left(\frac{\sqrt{x}}{\log x}\right) = (1.3727\ldots)\frac{\sqrt{x}}{\log x} + o\left(\frac{\sqrt{x}}{\log x}\right), \tag{34}$$

see [RG04, p. 7], [NW00, p. 342]. A lower bound for the counting function is given below.

**Lemma 14.** Let $x \geq 1$ be a large number, then $\#\{p = n^2 + 1 \leq x : p \text{ is prime}\} \gg x^{1/2}/\log x$.

*Proof*: By means of Theorem 1, the lower bound can be obtained by partial summation:

$$\#\{p = n^2 + 1 \leq x : p \text{ is prime}\} \gg \sum_{n^2+1 \leq x} \frac{\Lambda(n^2 + 1)}{n\sqrt{\log n}} \frac{n\sqrt{\log n}}{\log n} = \int_2^{x^{1/2}} \frac{t}{\sqrt{\log t}}\, d R(t) \gg \frac{x^{1/2}}{\log x}, \tag{35}$$

where $R(t) = \sum_{n^2+1 \leq t} \Lambda(n^2 + 1)\left(n\sqrt{\log n}\right)^{-1}$. ∎

## 9. Some Related Results

The problem investigated here can also be viewed as a special case $m = d$, $1 \leq d \leq 100$, of the results in [FI97 and [FI98] for primes of the forms $p = n^2 + m^2$, and $p = n^2 + m^4$, $m, n \in \mathbb{Z}$.

### 9.1 Counting Functions For Quadratic Primes In Two Variables

**Theorem 15.** ([FI98]) Let $f(r, s) = r^2 + s^4 \in \mathbb{Z}[r, s]$, (an absolutely irreducible polynomial over the integers), and let $x \geq 1$ be a large number. Then,

$$\sum_{n^2+m^4 \leq x} \Lambda(n^2 + m^4) = \frac{4\kappa}{\pi} x^{3/4} \left(1 + O\left(\frac{\log\log x}{\log x}\right)\right), \quad (36)$$

where the constant $\kappa = \int_0^1 \sqrt{1 - t^4} \, dt = \Gamma(1/4)^2 / (6\sqrt{2\pi})$.

**Theorem 15.** ([FI97]) Let $x \geq 1$ be a large number, and let $\chi \neq 1$ be a character mod 4. Then,

$$\sum_{n^2+m^2 \leq x} \Lambda(n) \Lambda(n^2 + m^2) = 2 \prod_{p \geq 2} \left(1 - \frac{\chi(p)}{(p-1)(p-\chi(p))}\right) x + O\left(\frac{x}{\log x}\right). \quad (37)$$

### 9.2 The PSI Function For Quadratic Polynomials

The well known psi function asymptotic formula

$$\psi(n) = \log \text{lcm}(1, 2, 3, \ldots, n) = n + O(n / \log n) \quad (38)$$

was extended to the quadratic arithmetic progression. The new result claims the following.

**Theorem 16.** ([CJ10[, [RZ13]) For $\theta < 4/9$, and $f(x) = x^2 + 1$, the expression

$$\psi_f(n) = \log \text{lcm}(1^2 + 1, 2^2 + 1, 3^2 + 1, \ldots, n^2 + 1) = n \log n + B n + O(n / \log^\theta n), \quad (39)$$

where the constant

$$B = \gamma - 1 - \frac{\log 2}{2} - \sum_{p \neq 2} \frac{\left(\frac{-1}{p}\right) \log p}{p - 1} = -.0662756342 \ldots. \quad (40)$$

### 9.3 Largest Prime Factors And Density Of Primitive Factors

Let $P(n)$ denotes the largest prime factor of the integer $n \geq 1$. Some of the history of the largest prime factor $P(n^2 + 1)$, and in general $P(f(n))$ for any irreducible polynomial $f(x) \in \mathbb{Z}[x]$, is discussed in [NW00, p. 345].

**Theorem 17.** ([ES90]) Let $f(x) \in \mathbb{Z}[x]$ be irreducible of degree $d = \deg(f) > 1$. Then there exists a constant $c_1 > 0$ such that for $x > x_1(f)$,

$$P\left(\prod_{n \leq x} f(n)\right) > x e^{e^{c_1 (\log\log x)^{1/3}}}. \tag{41}$$

**Theorem 18.** ([DI82]) For infinitely many integers $n \in \mathbb{Z}$, the largest prime factor of the polynomial $f(x) = x^2 + 1$ satisfies

$$P(n^2 + 1) \geq n^{1.202468\ldots}. \tag{42}$$

A factor $m \mid A_n$ of an integer sequence $A_n$, $n \geq 1$, is called a *primitive factor* if
(1) $m \nmid A_n$ for $n \leq n_0$, and (2) $m \mid A_n$ for some $n > n_0$.

Define the density function $\rho(f) = \lim_{x \to \infty} \#\{ n \leq x : f(n) \text{ has a primitive factor} \} / x$.

**Theorem 19.** ([EH08]) Let $x \geq 1$, and $f(x) = x^2 + d$ be an irreducible polynomial. Then, it satisfies

$$x \ll \rho(n^2 + d) \ll x. \tag{43}$$

### 10. Some Problems

**Problem 1.** Prove that there are infinitely many twin quadratic primes $n^2 + 1$, $n^2 + 3$ as $n \to \infty$. The sequence has the initial pairs (101, 103), (197, 199), (5477, 5479), (8837, 8839), ... . This problem is discussed in [NW00, p. 342].

**Problem 2.** Determine the minimal and maximal gaps of two consecutive quadratic primes $p_k = n^2 + 1$, $p_{k+1} = m^2 + 1$, $m, n \geq 1$. The minimal gap is $(n+2)^2 + 1 - (n^2 + 1) = 4n + 2 \geq \sqrt{p_k}$. And the average gap is also large:

$$x / \pi_f(x) = x / \left(c_0 x^{1/2} / \log x + o(x^{1/2} / \log x)\right) = c_0 x^{1/2} \log x + o(x^{1/2} \log x), \tag{44}$$

where $c_0 = 1.3727\ldots$ is a constant. Thus, this problem is quite different from the linear primes $\{2, 3, 5, 7, 11, 13, 17, 19, \ldots\}$, which has an average of $x / \pi(x) = \log x + o(\log x)$.

For $n^2 + 1 \leq 10\,000$, the sequence of primes is
2, 5, 17, 37, 101, 197, 257, 401, 577, 677, 1297, 1601, 2917, 3137, 4357, 5477, 7057, 8101, 8837, ...,

and the sequence of gaps is

3, 12, 20, 96, 60, 144, 176, 100, 620, 304, 1316, 220, 1220, 1120, 1580, 1044, 736, ... .

**Problem 3.** The exact values of several series associated with primes producing polynomials are known to be transcendental numbers:

$$\sum_{n\geq 0}\frac{1}{n^2+1}=\frac{\pi}{2}\frac{e^{2\pi}+1}{e^{2\pi}-1}+\frac{1}{2},\quad \sum_{n\geq 0}\frac{1}{n^2+3}=\frac{\pi}{2\sqrt{3}}\frac{e^{2\pi\sqrt{3}}+1}{e^{2\pi\sqrt{3}}-1}+\frac{1}{6},\quad \text{etc.,} \tag{45}$$

see [SR74, p. 189-199]. Does the transcendence of the series implies that the sequence of denominators contain infinitely many primes?

**Problem 4.** Prove that there are infinitely many primes $2p^2+1$ as the prime $p\to\infty$. The sequence $ap^2+1$, $a\geq p^\epsilon$, and $\epsilon>0$, has infinitely many primes as the prime $p\to\infty$, this sequence is constructed in [MK09].

**Problem 5.** Determine the least primitive root of the sequence of quadratic primes $n^2+1$ as $n\to\infty$. The primitive roots of some quadratic sequences are studied in [AS13].

**Problem 6.** Show that the function $f(s)=\sum_{n\geq 1}\Lambda(n^2+1)n^{-s}$ has a pole at $s=1$. Hint: Let

$$f(s)=\lim_{x\to\infty}\sum_{n\leq x}\frac{\Lambda(n^2+1)}{n^s}, \tag{46}$$

then use an approximation as in the proof of Theorem 1.

**Problem 7.** Prove that the class number $h(d)$ of the quadratic numbers field $\mathbb{Q}(\sqrt{d})$ satisfies $h(d)>1$ for infinitely many discriminants $d=n^2+1\geq 5$. This problem was considered in [MN88].